
\input amstex
\documentstyle{amsppt}
\topmatter
\title Differential-geometric methods for the lifting
problem and linear systems on plane curves \endtitle
\rightheadtext{DIFFERENTIAL-GEOMETRIC METHODS}
\author Emilia Mezzetti \endauthor
\affil Dipartimento di Scienze Matematiche - Universit\`a
di Trieste \endaffil
\address Piazzale Europa 1 - 34100
Trieste - Italy \endaddress
\email MEZZETTE\@UNIV.TRIESTE.IT (bitnet);
UTS882::MEZZETTE (decnet)\endemail

\abstract Let $X$ be an integral projective variety of
codimension two, degree $d$ and dimension $r$ and $Y$ be its
general hyperplane section. The problem of lifting
generators of minimal degree $\sigma$ from the homogeneous
ideal of $Y$ to the homogeneous ideal of $X$ is studied. A
conjecture is given in terms of $d$, $r$ and $\sigma$; it is
proved in the cases $r=1,2,3$. A description is given of
linear systems on smooth plane curves whose dimension is
almost maximal. \endabstract
 \endtopmatter
\document
\head Introduction \endhead
Let $X$ be an integral projective variety of dimension $r$ and
 degree $d$ in $\Bbb{P}^{r+2}$, the projective space of dimension
 $r+2$
over an algebraically closed field $K$ of characteristic $0$ and
 let
$Y=X\cap H$ be its general hyperplane section. We associate to $X$
the
integer $\sigma$, the minimal degree of a hypersurface in $H$
 containing
$Y$.

In this paper we
study the problem of lifting generators of degree $\sigma$
from the homogenous ideal of $Y$ to the homogenous ideal
of $X$. Answers in terms of relations between $d$ and
$\sigma$ are known for $r=1, 2$. For $r=1$, the Laudal's
\lq\lq generalized trisecant lemma'' (\lbrack L\rbrack, \lbrack
GP\rbrack) states that: if $X$ is an integral space curve of
 degree
$d>\sigma^{2}+1$, such that $Y$ is contained in a plane curve of
 degree
$\sigma$, then $X$ is contained in a surface of degree $\sigma$.
 For $r=2$,
there is the following analogous result (\lbrack MR \rbrack):
 if $X$ is an
integral linearly normal locally Cohen-Macaulay surface of
 $\Bbb{P}^{4}$,
of degree $d>\sigma^{2} - \sigma+2$, such that $Y$ is contained
 in a
surface of $\Bbb{P}^{3}$ of degree $\sigma$, then $X$ is contained
in a
hypersurface of degree $\sigma$.
In both cases, the bounds on the
degree are sharp and the border cases are classified
(\lbrack E\rbrack,
\lbrack St\rbrack, \lbrack M\rbrack).

In this paper we propose an approach to the problem, via a method of
differential geometry: our point of view, which is similar to that
of Gruson--Peskine (\lbrack GP\rbrack), consists in studying the
family of degenerate hypersurfaces of degree $\sigma$ containing the
hyperplane sections of $X$. In particular, we study the \lq\lq focal
linear system''
$\delta$ on a general member $G_Y$ of the family, a degenerate
hypersurface of degree $\sigma$ containing $Y=X\cap H$. The notion of foci of a
family of linear spaces or, more in
general, of varieties is classical; it was studied for example by Corrado
Segre in \lbrack Se\rbrack. Recently it was rediscovered and translated
in modern language in \lbrack CS\rbrack. Roughly speaking, the focal
linear system describes the intersections of $G_Y$ with the members of
the family which are \lq\lq infinitely near'' to it. In our case, $\delta$
is cut out on $G_Y$ by hypersurfaces of degree $\sigma +1$ containing $Y$,
out of $Y$. The important remark is the following: assume that $X$ is not
contained in a hypersurface of degree $\sigma$; then there is a lower
bound on the dimension of $\delta$. Now the strategy is to cut
successively with general hyperplanes, until we get  a linear system
of divisors of degree $\sigma (\sigma +1)-d$ and \lq\lq high''
dimension on an integral plane curve of degree $\sigma$.

The study of linear systems of maximal dimension on smooth plane
curves goes back to Max Noether; a modern proof of the classification of
such systems is given in \lbrack C\rbrack; in \lbrack H\rbrack \ similar
results are given for Gorenstein non-necessarily smooth curves. Here we
need informations about linear systems on integral plane curves whose
dimension is almost maximal with respect to the degree, such that the
general divisor is supported in the smooth locus. Since there are no exhaustive
results on these systems in the literature,
it has been necessary to devote part of the paper to this subject. In
particular, we prove a theorem giving a complete description of linear
systems whose dimension is one less than the maximal; it seems to us of
independent interest and susceptible of further extensions. We plan to
come back to this subject in a successive paper.

If the dimension of $X$ is $r\leq 3$, thank to this result, we are able
to find an upper bound, depending on $\sigma$, on the degree $d$ of $X$.
Precisely, for curves and surfaces we refind the known results. For
$3-$folds, we prove the following (Theorem 4.10):

\proclaim\nofrills{} Let $X$ be an integral non--degenerate subvariety of
$\Bbb P^5$ of dimension $3$ and degree $d$. Let $\sigma >5$ be the minimal
degree of a degenerate hypersurface containing a general hyperplane
 section of $X$. If $d>\sigma ^2-2\sigma +4$, then $X$ is contained in a
hypersurface of degree $\sigma$.\endproclaim
Now the question arises of finding a reasonable conjecture for the
general case. We propose the following:

\proclaim{Conjecture} Let $X$ be an integral projective variety
of dimension $r$ and degree $d$ in $\Bbb{P}^{r+2}$, $Y=X\cap H$ be its
general hyperplane section and $\sigma$ be the minimal degree of a
hypersurface in $H$ containing $Y$. If $X$ is not contained in a
hypersurface of degree $\sigma$, then $d\leq \sigma^{2} - (r - 1)\sigma +
\binom r2 +1$. \endproclaim

Note that, if $X$ and $\sigma$ are defined as above and
$h^0(\Cal I _X(\sigma ))=(0)$, according to the construction of
Gruson--Peskine (\lbrack GP\rbrack) there is an exact
sequence of sheaves on a general hyperplane $H$, of the form $$ 0 \to N
\to \Omega_{H}(1) \to J(\sigma) \to 0,$$  where: $N$ is a reflexive sheaf
of rank $r+1$, $\Omega_{H}$ is the cotangent bundle of $H$ and $J $ is an
ideal sheaf of $\Cal O_G$, where $G$ is a hypersurface of $H$ containing
$Y$; $J$ defines a closed subscheme $\Delta$ of $G$, containing $Y$, of
dimension $r - 1$ and degree $\delta \ge d$.

By computing the Chern classes, we get: $c_2(N(1)) = \sigma^{2} -
(r - 1)\sigma + \binom r2 +1 - \delta$,  so the conjecture would follow
from the condition $c_2(N(1))\geq 0$.

 Note that if
this conjecture were true, then the bound on the degree would be sharp.
In fact there is a whole class of examples of integral varieties of
dimension $r$ and degree   $d=\sigma^{2} - (r-1)\sigma + \binom r2 +1$,
not lying on a hypersurface of degree $\sigma$, generalizing the examples
in dimension $1$ and $2$. They are arithmetically Buchsbaum varieties
with $\Omega$ - resolution: $$0 \to r\Cal O_{\Bbb{P}^{r+2}} (-1)\oplus
\Cal O_{\Bbb{P}^{r+2}}(r-1-\sigma) \to \Omega_{\Bbb{P}^{r+2}}(1) \to \Cal
I_X(\sigma) \to 0$$  (see \lbrack Ch \rbrack ). Note that in these
examples the only non-zero intermediate cohomology for $\Cal I_X$ is
$h^1(\Cal I_X(\sigma - 1))=1$, i.e. the situation is \lq\lq the best
possible\rq\rq.

We believe that our approach to the problem, if suitably deepened, can
be useful to study this conjecture also for varieties of higher
dimension.

The paper is organized as follows: in \S 1 we introduce
 the concept of foci and focal divisors for families of degenerate
hypersurfaces containing the hyperplane sections of $X$. In \S 2 we
introduce the focal linear system $\delta$ on a general hypersurface
$G_Y$ of degree $\sigma$ containing $Y=X\cap H$. We prove that, if
$h^0(\Cal I_X(\sigma ))=0$ then $dim \ \delta \geq 2r+h^0(\Cal
I_Y(\sigma ))-3$; if moreover $\sigma > \frac {1}{r+1} \binom {2r}r$ ,
then $dim \ \delta \geq 2r+h^0(\Cal I_Y(\sigma ))-2$. In \S 3
we give a complete description of the linear systems on smooth
plane curves, whose dimension is either the maximal with respect
to the degree or the maximal minus one (Theorem $3.2$). Thanks to this
result, in \S 4 we are able to prove the above conjecture for $r=1, 2$ and
 for $r=3$ under the assumption $\sigma >5$.

\remark{Aknowledgements}  This work was
 done in the framework of the research project \lq\lq Lifting problem,
Halphen problem and related topics\rq\rq of Europroj. I am grateful to all
participants to this project for stimulating talks. In particular, I
would like to thank Luca Chiantini and Ciro Ciliberto for suggesting the
use of focal linear systems. \endremark

\head 1. Families of degenerate hypersurfaces
\endhead
We work in $\Bbb{P}^{r+2}$, $r\ge 1$, the projective space of
dimension $r+2$ over an algebraically closed field $K$ of characteristic
$0$. In the following it will be denoted simply by $\Bbb{P}$. By
 $\check \Bbb P$ we will denote
 the dual projective space; if $L\subset\Bbb P$ is a linear subspace,
$\check L$ will denote its dual, i.e. the subspace of $\check\Bbb P$ of
hyperplanes containing $L$.

Let $X\subset \Bbb{P}$ be an
integral non-degenerate projective variety
of dimension $r$. Let $H$ be a general hyperplane and $Y= X\cap H$ be a
general hyperplane section of $X$. Let $\Cal I_Y\subset\Cal O_H$ be the
ideal sheaf of $Y$ in $H$. We associate to $X$ the integers $\sigma : =
min \lbrace k\in \Bbb Z \mid h^0( \Cal I _Y (k)) \ne 0 \rbrace \ge 2$ and
$n : = h^0(\Cal I _Y (\sigma ))\ge 1$. Note that every element $G$ in
$H^0(\Cal I _Y(\sigma ))$ is integral, by the assumption \lq\lq $X$
integral\rq\rq  and the minimality of $\sigma$; moreover, $Y$ is not
contained in the singular locus of $G$.

If $H$ lies   in a suitable non--empty
open subset $U$ of $\check \Bbb P$, the  hypersurfaces
 of $H$ of degree $\sigma$  containing $X\cap H$ form a linear system of
dimension $n-1$. So the family of hypersurfaces of $H$ of degree
$\sigma$ containing  $X\cap H$, for $H$ varying in $U$, is
parametrized by an irreducible variety $T$ of dimension $(r+2) + (n-1) =
r+n+1$. If $t\in T$, let $G_t $ be the corresponding hypersurface (we
denote in the same way
 a hypersurface and
a polynomial defining it); we denote by $\Sigma _T$ the
following incidence family $\Sigma _T \subset T \times \Bbb P$ and
call it the total family:   $$ \CD \Sigma _T = \lbrace (t,P)| P \in G_t
\rbrace   @>p_1>> T  \\ @VVp_2V   @.   \\ \Bbb P   @. \endCD $$
 where $p_1$ and $p_2$ are the projections.

In the following we will study some subfamilies $\Sigma _Z$
of $\Sigma_T$ constructed as follows: if  $h^0(\Cal I_Y(\sigma))=n=1$, let
$G_Y$ be the unique hypersurface of degree $\sigma$ in $H$ containing
$Y=X\cap H$; if $n >1$, let us fix $l_1,...,l_{n -
1}$, $n-1$ general lines of $\Bbb P$; then $l_i\cap H = P_i$,
$i=1,...,n-1$, are general points in $H$: we call $G_Y$ the unique
hypersurface of degree $\sigma$ in $H$ containing $Y$ and passing through
$P_1,...,P_{n-1}$. Let $Z\subset \check \Bbb P$ be a smooth integral
closed subvariety; we define the incidence family $\Sigma _Z \subset Z
\times \Bbb P$ :
 $$
\CD
\Sigma _Z = \lbrace (H,P)| H \in Z, P \in G_Y, Y=X\cap H \rbrace  @>p_1>> Z  \\
@VVp_2V   @. \\
\Bbb P   @.
\endCD
 $$
 From now on we suppose that, if $H$ is general
in $Z$, $G_Y$ is uniquely determined and irreducible, non--degenerate;
for such $H$,
 $dim \ p_1^{ - 1}(H) = r$ hence $dim \ \Sigma _Z =
dim \ Z + r$ and $dim \ p_2(\Sigma _Z) \le dim \ Z + r$.
\proclaim{1.1.
Proposition} Let  $Z= \check \Bbb P^h \subset \check \Bbb P$ be a linear
subspace of dimension at least two, $\Sigma_Z$ the incidence family
constructed as above. If  the projection $p_2: \Sigma_Z \to \Bbb P$ is
not dominant, then $dim \ p_2(\Sigma_Z) = r+1$ and $deg \ p_2(\Sigma_Z) =
\sigma$. \endproclaim
\demo{Proof} Choose $H$ general in $Z$: we prove
that $p_2(\Sigma_Z)\cap H = G_Y$, $Y=X\cap H$. It is obvious that
$G_Y\subset p_2(\Sigma_Z)\cap H$. Let $P$ be a general point of
$p_2(\Sigma_Z)\cap H$: since $p_2$ is not dominant, $dim \ p_2^{-1}(P) >
dim \ Z + dim \ p_1^{ - 1} (H) - r - 2 = dim \ Z - 2$; we have: $p_2^{ -
1}(P) \simeq p_1(p_2^{ - 1}(P)) \subset \check P\cap Z$. Moreover $dim
\ \check P\cap Z = dim \ Z -1$, because
$P\notin \Bbb P^h$ (if $P\in\Bbb P^h$, then
$ p_2(\Sigma_Z)\cap H\subset\Bbb P^h$ for
any $H$ in $Z$, so  $p_2(\Sigma_Z)\cap
H=\Bbb P^h$ and $G_Y\subset\Bbb P^h$; this
implies $r\leq h$, which contradicts $dim \
Z=r-h-1\geq 2$).
                   So $p_2^{ - 1}(P) = \check P\cap Z$, hence $H\in
p_1(p_2^{ - 1}(P))$ and $P\in G_Y$. Therefore $p_2(\Sigma_Z)\cap H
\subset G_Y$. Finally note that the intersection $p_2(\Sigma_Z)\cap H$ is
reduced: otherwise $H$ would be a general tangent hyperplane to
$p_2(\Sigma_Z)$ and $G_Y$ the tangency locus, which should be linear.
\enddemo
\proclaim{1.2. Corollary} Let $Z\subset \check \Bbb P$  be as in
$1.1$.  If $H^0(\Cal I_X(\sigma))=(0)$, the map $p_2$ is dominant.
\endproclaim

{}From now on we consider families $\Sigma _Z$ such that $p_2$ is dominant.

Following the exposition  in \lbrack CS\rbrack, we will introduce
now the concept of foci for a family $Z$ constructed as above or
for the total family. We denote
 $T(p_2) = Hom (\Omega_{Z\times \Bbb P \mid \Bbb P}, \Cal O_{ Z\times \Bbb
P})$ and call {\it global characteristic map} of the family $\Sigma_Z $
the map $d: T(p_2)\vert _ {\Sigma_Z}\to \Cal N_{\Sigma_Z, Z\times \Bbb
P}$ which is defined by the commutative diagram:
$$ \CD {} @.   {}   @.
0  @.  {}  @.  \\
 @.      @.   @VVV  @.  @.  \\
{}  @.   {}  @.  T(p_2)\vert _ {\Sigma_Z}  @>d>>
 \Cal N_{\Sigma_Z , Z\times \Bbb P}  @>>>  0  \\
@.      @.   @VVV  @\vert    @.  \\
0  @>>>  T_{\Sigma _Z}   @>>>  T_{Z\times \Bbb P \vert \Sigma_Z}
 @>>> \Cal N_{\Sigma_Z , Z\times \Bbb P}  @>>>  0  \\
@.   @VVV   @VVV   @.   @.   \\
{}   @.   p_2^*T_{\Bbb P \vert\Sigma_Z}  @=  p_2^*T_{\Bbb P \vert\Sigma_Z}  @.
{}   @.      \\
\endCD
\tag 1.3 $$

For any $H\in Z$, let us consider the restriction of $d$ to the
 fiber $p_1^{-1}(z)\simeq G_Y$, $Y=X\cap H$:
$$
\CD
d_Y: T_{H,Z}\otimes \Cal O_{G_Y}   @>>>   \Cal N_{G_Y,\Bbb P}  \\
 				@\vert    @\vert  \\
	  \Cal O_{G_Y}^{dim Z}   @>>>   \Cal O_{G_Y}(\sigma)\oplus \Cal O_{G_Y}(1)
 \endCD
$$
Note that the map $d_Y$  is induced by the characteristic map
 which associates to an infinitesimal deformation of $G_Y$ in the family
$\Sigma_Z$ parametrized by $Z$ the corresponding deformation of $G_Y$ in
the Hilbert scheme of $\Bbb P^{r+2}$. Analogous notations are used when
considering the total family $\Sigma _T$.

\example{1.4. Example}  Let
$H\in Z$ be the hyperplane of equation $x_{r+2}=0$. The map
$d_Y$ is given by a $(2\times dim \ Z)$-matrix $M$. If we consider the
total family $\Sigma_T$, we may suppose that $M$ is of the form $M =
\left( \matrix  F_0  & F_1 &  \hdots & F_{r+1} & F_{r+2} & \hdots  &
F_{r+n} \\
 x_0 & x_1 & \hdots & x_{r+1} & 0 & \hdots & 0  \endmatrix \right),$
where $F_i \in
H^0(\Cal O_{G_Y}(\sigma ))$, $i=0,...,r+n$. If we fix lines
$l_1,...,l_{n-1}$ as before and take $Z=\check \Bbb P^{r+2}$ (resp.
$Z=\check P, P=(0:\hdots :0:1:0))$ we may suppose that $M$ is of the form
$M= \left( \matrix  F_0  & F_1 &  \hdots & F_{r+1}  \\
x_0 & x_1 & \hdots & x_{r+1}
 \endmatrix
\right)$ (resp. $M= \left( \matrix
F_0  & F_1 &  \hdots & F_{r}\\
x_0 & x_1 & \hdots & x_{r}    \endmatrix \right)$ ).\endexample

\definition{1.5. Definition} Let $\Sigma_Z$ denote the total family or a
subfamily constructed in the usual way. Let $d_Y$ be the
restriction of the global characteristic map to $G_Y$, $Y=X\cap H$. The
points of the support of $coker \ d_Y$ are called the {\it foci of $G_Y$
in the family $\Sigma_Z$}. We shall denote them by $F_Y$. A matrix $M$
associated to the map $d_Y$ is called a {\it focal matrix} of $G_Y$ in
$\Sigma_Z$.\enddefinition

$F_Y$ is the closed subset of $G_Y$ defined by the vanishing of the
maximal order minors of $M$. In the example $1.4$ above, in the cases
$Z=\check \Bbb P^{r+2}$  or $Z=\check P$, $F_Y$ is formed by the points
of $G_Y$ where  all the
polynomials $F_ix_j - F_jx_i$, $i\ne j$, are zero.

If $H\in Z$, for any linear subspace $\pi$ of codimension $2$ in $H$,
the inclusion $$T_{H,\check \pi} \hookrightarrow T_{H,Z}$$  induces a
diagram which defines a map $d_{Y,\pi}$:
$$
\CD
0  @.  {}  \\
@VVV   @.   \\
\Cal O_{G_Y}^2   @=  \Cal O_{G_Y}^2 \\
@VVV            @VVd_{Y,\pi}V \\
\Cal O_{G_Y}^{dimZ}   @>>d_Y>   \Cal O_{G_Y}(\sigma)\oplus \Cal
O_{G_Y}(1). \endCD
$$
\definition{1.6 Definition} The points of the support of $coker \
d_{Y,\pi}$  are called the {\it foci of $\pi$ on $G_Y$ in the family
$\Sigma_Z$}. We shall denote them by $F_{Y,\pi}$. They either form a
divisor on $G_Y$, which is called the {\it focal divisor of $\pi$}, or
fill $G_Y$; in this case we say that the focal divisor is undetermined.
\enddefinition

Obviously we have $F_Y\subset F_{Y,\pi}$ for any
$\pi\subset H$ and $F_Y=\cap F_{Y,\pi}$.

\proclaim{1.7. Proposition} In
the above situation:
\roster \item $P\in F_{Y,\pi}$ if and only if the
linear space $<P, \pi >\check {}$  is tangent at $H$ to the image via
$p_1$ of the fiber $ p_2^{-1}(P)$;
\item $P\in F_Y$ if
and only if $(H,P)$ is a singular point of $p_2^{-1}(P)$ or  $dim \
p_2^{-1}(P)> dim \ Z-2$. \endroster \endproclaim
\demo{Proof}
(1)  $P\in F_{Y,\pi}$ if and only
 if the morphism $d_{Y,\pi}(P)$ induced by $d_{Y,\pi}$ on the stalks at
$P$ is not isomorphic, i.e. there exists a non-zero $\tau \in
T_{H,\check\pi}$  such that the section $d_{Y,\pi}(\tau\otimes\Cal
O_{G_Y})$ vanishes at $P$. Then by the diagram (1.3) we have the claim.

(2) By (1), $P\in F_Y$ if and only if, for $\pi$ general in $H$,
$p^{-1}(< P, \pi>)\check  {}$ is tangent to the fiber $p_2^{-1}(P)$ at
$(H,P)$.  \enddemo

\proclaim{1.8. Corollary} If $p_2$ is dominant, then
the inclusion $F_Y\subset G_Y$ is strict. Moreover $Y\cup P_1\cup ...\cup
P_{n - 1} \subset F_Y$.  \endproclaim
\demo{Proof} The first assertion
follows because the general fiber of $p_2$ is reduced; the second
one because $dim \ p_2^{-1}(Q) = dim \ Z - 1$ if $Q\in Y$ or $Q=P_i,
i=1,...n - 1$. \enddemo

\head 2. The focal linear system  \endhead
In this section we will study the foci of the families $\Sigma _Z$ of \S
1.

Let $X$ be an integral non--degenerate variety of codimension $2$ in
$\Bbb P^{r+2}$. Let $\Sigma _Z$ denote the total family of degenerate
hypersurfaces of minimal degree $\sigma$ containing the hyperplane
sections of $X$, or a subfamily of its as in \S 1. Let $H$ be a
general hyperplane, $Y=X\cap H$; let $d_Y$ be the restriction of the
global characteristic map of $\Sigma _Z$ to $G_Y$:
$$d_Y: \Cal O_{G_Y}^{dim Z} \longrightarrow \Cal N _{G_Y,\Bbb
P}\simeq \Cal O _{G_Y}(\sigma)\oplus\Cal O_{G_Y}(1).$$
Let us consider the map $\wedge ^2d_Y: \wedge ^2\Cal O_{G_Y}^{dim
\ Z}\rightarrow det \ \Cal N_{G_Y,\Bbb P}\simeq\Cal O_{G_Y}(\sigma +1)$.

\definition{2.1. Definition} The {\it focal linear system} of the
family $\Sigma _Z$ on $G_Y$ is the projectivized of the image of
$\wedge ^2d_Y$. We will denote it by $\delta _Z$ (or simply by $\delta$).
Note that, by 1.8, $\delta$ is cut out on $G_Y$ by hypersurfaces of
$H$ of degree $\sigma +1$ containing $Y$, i.e. $\delta\subset\Bbb
P(\Cal O_{G_Y}(\sigma +1)(-Y))$.\enddefinition

Let us fix general lines $l_1,\dots ,l_{n-1}$ and consider the
corresponding family $\Sigma _Z$, $Z=\check\Bbb P^{r+2}$. Let $H$ have
equation $x_{r+2}=0$ and let $V=H^0(\Cal O_H(1))$ be the $K$-vector
space with basis $x_0, ... ,x_{r+1}$.  Consider $M$, a focal matrix
of $G_Y$ in $\Sigma _Z$. Let $\phi_M :
\Lambda ^rV^*  \to  H^0(\Cal O_{G_Y}(\sigma +1))$ be the linear map which
acts in the following way: for a permutation
$x_{i_0}^*,...,x_{i_{r+1}}^*$ of the elements of the dual of the fixed
basis of $V$, $\phi_M(x_{i_0}^*\wedge ... \wedge x_{i_{r-1}}^*) $ is  the
minor of $M$ corresponding to the columns of indices $i_r, i_{r+1}$.
Then $\delta _Z=\Bbb P(Im \ \phi _M)$.

Let $\Bbb G = \Bbb G(r-1,r+1)$ be the grassmannian of $(r-1)$-planes in
$H$; it is a projective variety of dimension $2r$ embedded in
$\Bbb P(\wedge ^rV^*)$ via the Pl\"ucker embedding. So $\phi _M$ restricts
to a rational map from $\Bbb G$ to the focal linear system $\delta$; by
definition, $\phi _M(\pi)=F_{Y,\pi}$ (if it exists).

\proclaim{2.2. Lemma.} Let $V$ be a hypersurface of $\Bbb P^n$; if $V$
contains  a family of dimension $2$ of
$(n-2)$--planes, then it is a hyperplane.\endproclaim
\demo{Proof} The intersection of $V$ with a general
linear space of dimension $3$ is a surface of $\Bbb P^3$ containing a
family of lines of dimension $2$, hence a plane.\enddemo

\proclaim{2.3. Theorem}Let $X$ be an integral non--degenerate variety of
codimension $2$ in $\Bbb P^{r+2}$, $H$ a general hyperplane, $Y=X\cap
H$. Let $\sigma$ be the minimal degree of a degenerate hypersurface
containing $Y$; suppose that $H^0(\Cal I_X(\sigma ))=(0)$. Let
$Z=\check\Bbb P^{r+2}$ and $\Sigma :=\Sigma _Z$ be
the incidence family constructed in \S 1. Let
$\delta$ be the focal linear system on $G_Y$; then $dim \ \delta
\geq 2r-2$. If  moreover $\sigma >deg \ \Bbb G =  \frac {1}{r+1}
\binom {2r}r$, then $dim \ \delta \geq 2r-1$.\endproclaim
\demo{Proof} We have  to show that,
fixed $2r-2$ arbitrary points of $G_Y$, there is an element of $\delta$
through them.

Note that, if we fix a point $P$ in $G_Y$, the linear system $\delta _P:
=$ \lbrace divisors of $\delta$ passing through $P$\rbrace comes via
$\phi _M$ from a hyperplane section of $\Bbb G$. So $2r-2$ points
of $G_Y$ determine the intersection  $W$ of $\Bbb G$ with a linear space
of codimension $2r-2$; since $dim \ \Bbb G=2r$, $W$ has dimension at least
$2$.  It is enough to show that the focal divisor associated to $\pi$
cannot be undetermined for each $\pi$ in $W$.

The union of the $(r-1)$-planes
of $W$ is $H$ or a hyperplane in $H$, by  lemma 2.2. We distinguish
now two cases:

a) the $(r-1)$-planes of $W$ are two by two in general
position.

Assume $H$ has equation $x_{r+2}=0$. We choose $\lbrack (r+2)/2\rbrack$ of
the $(r-1)$-planes of $W$: $\pi _0, \pi _1, \dots$, not contained in
$G_Y$; we may suppose that they have equations $\pi _0: x_0=x_1=0, \pi _1:
x_2=x_3=0,\dots , \pi _i: x_{2i}=x_{2i+1}=0,\dots$. The focal divisor
associated to $\pi _i$ is given by the equation
$x_{2i}F_{2i+1}-x_{2i+1}F_{2i}=0$; we assume by contradiction that
it is undetermined for each $i$. Since it is undetermined for $i=0$, then
$x_1F_0-x_0F_1=LG_Y$, $L$ a linear form; since $G_Y\notin (x_0, x_1)$, we
get $L\in(x_0, x_1)$ so $x_1F_0-x_0F_1=(ax_0+bx_1)G_Y$  $(a,b\in K)$,
i.e. $x_1(F_0-bG_Y)=x_0(F_1+aG_Y)$. Hence there exists a polynomial $A$
such that $F_0=x_0A+bG_Y$, $F_1=x_1A-aG_Y$; so in the
matrix $M$ we may replace $F_0$ by $x_0A$ and $F_1$ by $x_1A$.  Similarly
we may assume $F_2=x_2A'$, $F_3=x_3A'$, etc. The minors $x_0x_3(A-A')$,
$x_0x_2(A-A')$, ...  define hypersurfaces containing $Y=X\cap H$; by the
minimality of $\sigma$, we conclude that $M$ has the following form:

$ \left( \matrix  x_0A  &   \hdots & x_{r+1}A  \\
x_0 & \hdots & x_{r+1}
 \endmatrix
\right)$ if $r$ is
even;

$\left( \matrix  x_0A  &  \hdots & x_rA  & F_{r+1}  \\
x_0  & \hdots & x_r  & x_{r+1}
 \endmatrix
\right)$ if $r$ is odd.

In the case \lq\lq $r$ even'', we have a contradiction, because by 1.2
$p_2$ is dominant and by 1.8 $F_Y$ is a proper subset of $G_Y$. In the
case \lq\lq $r$ odd'', let us remark that the matrix
$M':= \left( \matrix  x_0A  &   \hdots & x_{r}A  \\
x_0 & \hdots & x_{r}
 \endmatrix
\right)$ is a focal matrix for
the family $\Sigma _{Z'}$, $Z'=\check P$, ($P$ the point of coordinates
$(0:\dots :0:1:0)$); since $rk \ M'<2$, again by 1.2 and 1.8 we have a
contradiction.

b) the $(r-1)$-planes of $W$ are all contained in a $r$-plane or
intersect each other along a fixed $(r-2)$-plane $L$.

Assume by contradiction that the focal divisor is undetermined for any
$\pi$ in $W$. In the first case arguing as in a) we find that $M$ has the
form $ \left( \matrix  x_0A  &   \hdots & x_{r+1}A  \\
x_0 & \hdots & x_{r+1}
 \endmatrix
\right)$  and conclude in the
same way. In the second case, we may assume that $L$ has equations
$x_0=x_1=x_2=0$ and choose two $(r-1)$-planes of equations $x_0=x_1=0$,
$x_0=x_2=0$; we find that $M$ has the form: $ \left( \matrix  x_0A  &
x_1A &x_2A & F_3 &
 \hdots & F_{r+1}  \\
x_0 & x_1 & x_2 & x_3 & \hdots & x_{r+1}
 \endmatrix
\right)$. Let $Z'=\check L$: then
$ \left( \matrix  x_0A  &   x_1A & x_2A  \\
x_0 & x_1 & x_2
 \endmatrix
\right)$ is a focal matrix for $\Sigma _{Z'}$,
which contradicts 1.2 and 1.8.

Assume now that  $\sigma >deg \ \Bbb G$ and  fix $2r-1$ arbitrary points
of $G_Y$; they determine the intersection  $W$ of $\Bbb G$ with a linear
space of codimension $2r-1$; if the intersection is proper, then $W$ is a
curve of degree $deg \ \Bbb G$ parametrizing a family of dimension $1$ of
$(r-1)$-planes of $H$; otherwise we have a family of bigger dimension. In
the first case, the union of the $(r-1)$-planes of $W$ is a hypersurface
$W'$ in $H$ of degree $deg \ \Bbb G$; since $\sigma >deg \ \Bbb G$ and
$G_Y$ is irreducible,$ W'\not\subset G_Y$. In the second case we find
again $W=H$ or a hyperplane of $H$. The discussion of the two possible
cases goes exactly as above. \enddemo

\proclaim{2.4. Corollary} Let $X$ be an integral non--degenerate variety
of codimension $2$ in $\Bbb P^{r+2}$, $H$ a general hyperplane, $Y=X\cap
H$. Let $\sigma$ be the minimal degree of a degenerate hypersurface
containing $Y$; suppose that $H^0(\Cal I_X(\sigma ))=(0)$. Let $\delta
_T$ be the focal linear system of the total family on $G_Y$.  Then
$dim \ \delta _T\geq  32r+n-3$; moreover, if $\sigma >\frac {1}{r+1}
\binom {2r}r$,
 $dim \ \delta _T \geq 2r+n-2$.\endproclaim
\demo{Proof} Observe that when we
choose $n-1$ general lines, as made above to construct the family
$\Sigma _Z$ of theorem 2.3,  we in fact impose to the general divisor of
$\delta _T$ the passage through $n-1$ general points.\enddemo

\head 3. Special linear systems on integral plane curves \endhead

In this section we will study linear systems on plane curves of
\lq\lq almost maximal'' dimension.

 Let $C$ be an integral plane curve of
degree $d>3$; let $g^r_n$ be a complete linear series on $C$ such that
the support of its general divisor $D$ is contained in the smooth locus of
$C$, so $g^r_n$ corresponds to an invertible sheaf on $C$. Let $H$ be a
linear section of $C$, $K \equiv (d-3)H$ the divisor associated to the
dualizing sheaf. Suppose that $g^r_n$ is special. We may write
$n=\alpha d-\beta$, $\alpha\leq d-3$, $0\leq \beta <d$. Then:
\proclaim{3.1. Theorem
(see \cite C, \cite H)}
In the above situation,  $r=r(n,d)-\varepsilon$, $\varepsilon \geq 0$,
where
$$r(n,d)=\cases \frac{1}{2}\alpha (\alpha +3)-\beta &\text{ if }
\alpha\geq\beta -1 \\
\frac{1}{2}(\alpha -1)(\alpha +2) &\text{ if } \alpha\leq\beta
-1.\endcases$$ \endproclaim

In the quoted papers, the linear systems with
$r=r(n,d)$ are completely described. Next theorem extends these results,
giving the description of linear systems with $\varepsilon\leq 1$.
\proclaim{3.2.
Theorem} If $\varepsilon\leq 1$ and $\varepsilon +1\leq \alpha\leq
 d-\varepsilon -3$, then
one of the following happens:
\roster
\item"(i)" $g^r_n 	=\mid \alpha H-D+F\mid$, $\alpha \geq \beta
-1-\varepsilon$, where $D$ and $F$ are effective disjoint divisors and
	{%
  \itemitem{(i1)} $deg \ D = \beta +\varepsilon$, $deg \ F=\varepsilon$ if
$\alpha \geq \beta -1$,
  \itemitem{(i2)}  $deg \ D = \beta$, $deg \ F=0$ if $\alpha =\beta -2$,
$\varepsilon=1$;
  \par}
\item"(ii)" $g^r_n 	=\mid (\alpha -1)H-D+F\mid$, $\alpha \leq \beta
-1+\varepsilon$, where $D$ and $F$ are effective disjoint divisors and
	{%
  \itemitem{(ii1)} $deg \ D =\varepsilon$, $deg \ F = d-\beta
+\varepsilon$ if $\alpha \leq \beta -1$,
  \itemitem{(ii2)}  $deg \ D = 0$, $deg J\ F =d-\beta$ if $\alpha
=\beta$, $\varepsilon =1$;   \par}
\item"(iii)" $g^r_n 	=\mid (\alpha +1)H-P_1-\dots -P_9\mid$, if
$\varepsilon =1$, $\alpha =2$, $\beta =3$, $d=6$, and $P_1,\dots ,P_9$ are
the complete intersection of two cubics. \endroster
 $F$ is a
fixed divisor; in (i) $D$ imposes independent conditions to the curves of
degree $\alpha$              .\endproclaim

\remark{3.3. Remark}  Note that the assumption on $D$ of imposing
independent conditions to the curves of degree $\alpha$ only means that
$D$ is not contained in a line in the cases i1) when $\varepsilon =1$,
$\alpha =\beta -1$, and i2). One can easily see that dually the condition
\lq\lq $F$ fixed'' means that $F$ is not contained in a line in the cases
ii1) when $\alpha =\beta -1$ and ii2).\endremark

 The proof of the theorem
is based on some lemmas.

If $\Cal L$ is a linear series on $C$ and $A$ an
effective divisor, $c(A,\Cal L$) will denote the number of independent
conditions imposed by $A$ to the divisors of $\Cal L$ and $\Cal L(A)$ the
linear series generated by the divisors of the form $D+A$, where
$D\in\Cal L$.
\proclaim{3.4. Lemma} Let $A$ be an effective divisor on $C$, $\Cal L$ a
linear series. Then $A$ is not fixed for the series $\Cal L(A)$ if and
only if $c(A,K-\Cal L) < deg \ A$; moreover: $dim \ \Cal L(A)-dim \ \Cal
L = deg \ A-c(A,K-\Cal L)$.\endproclaim
 \demo{Proof} It
is just an application of the theorem of Riemann-Roch.\enddemo
\proclaim{3.5. Lemma}
Let $E$ be a group of $m$ points in the plane; define $\tau
= \ max\lbrace t\in\Bbb Z\mid h_E(t)<m\rbrace$ where $h_E$ is the Hilbert
function of $E$. If $s$ is an integer such that $m\geq s^2$ and
$\tau\geq \frac{s^2-3s+m}{s}$, then one of the following happens:
\roster
\item"(i)" there
is a $t$, $0<t<s$, and a subset $E'$ of $E$ contained in a curve of degree
$t$, with $$t(\tau -t+3)  \leq deg \ E' \leq t(\tau +\frac{5-t}2);$$
\item"(ii)" $\tau =\frac{s^2-3s+m}s$  and $E$ is a
complete intersection of type $(s, m/s)$. \endroster\endproclaim
\demo{Proof} See \cite {EP}.\enddemo
\remark{3.6.Remarks}

 1) If $\alpha =1$, i.e. $n=d-\beta$,
we may have $\varepsilon >0$ only if $\beta =0,1$. If $\beta =1$, then
$r(n,d)=1$, so if $\varepsilon =1$, the series is a $g^0_{d-1}$ whose
points are not all on a line; if $d-2$ of these points are on a line, then
we may write $g^0_{d-1}=\mid H-P_1-P_2+Q\mid$. If $\beta =0$, then
$r(n,d)=2$. If $\varepsilon =1$, the series is a $g^1_{d}$: if it has a
fixed point $P$, then $g^1_{d}(-P)=g^1_{d-1}=\mid H-Q\mid$, so
$g^1_d=\mid H-Q+P\mid$; if it has no fixed point, we may assume that $D$
is a divisor of the series formed by distinct points; we have: $dim
\mid K\mid -dim \mid K-D\mid = d-1$, i.e. $D$ imposes $d-1$ conditions
to the curves of degree $d-3$; by 3.5 $D$ contains $d-1$ points lying on a
line, so $\mid D\mid = \mid H-Q+P\mid$.  If $\varepsilon =2$, the series
is a fixed divisor $D$ of degree $d$, which imposes independent conditions
to $\mid K\mid $, so it does not contain $d-1$ points on a line.

2) A
special $g^r_n$ satisfies the condition \lq\lq $r>n-p_a(C)$'' (where
$p_a(C)$ denotes the arithmetic genus of $C$).  If
$r=r(n,d)-\varepsilon$, this is equivalent to the condition
$$\varepsilon \leq\frac{d^2-(2\alpha +3)d+\alpha(\alpha +3)}2$$ in the
case $\alpha \geq \beta -1$, and to the condition
$$\varepsilon\leq\frac{d^2-(2\alpha +3)d+\alpha ^2+\alpha -2}2 +\beta$$ in
the case $\alpha \leq \beta -1$. In particular, if $\alpha =d-3\geq
\beta -1$, then always $\varepsilon =0$; if $\alpha =d-3$, $\beta =d-1$,
then $\varepsilon =0$ or $1$.  If $\varepsilon =1$ and $D$ is a general
divisor of the series, then $deg \mid K-D\mid = d-1$, $dim \mid K-D\mid
=0$. So $\mid K-D\mid =F$, a fixed divisor of degree $d-1$, and $\mid
D\mid = \mid (d-3)H-F\mid$; if moreover $F$ contains $d-2$ points on a
line, then $\mid D\mid = \mid (d-4)H-Q+P_1+P_2\mid$.\endremark

\smallpagebreak
In the following the symbol $\oplus$ will denote the
operation of minimal sum of linear series.
 \proclaim{3.7. Lemma} With the above
notations, let $g^r_n$ be a special linear series on $C$ with
$r=r(n,d)-\varepsilon$, $\varepsilon >0$. \roster
\item"a)"
If $\alpha\geq\beta -1$, $\alpha\leq d-\varepsilon -3$, then
$c(H,g^r_n)\leq \alpha  +1$;
 \item"b)" if $\alpha  \leq \beta   -2-\varepsilon$    , then
$c(H,g^r_n)\leq \alpha$;
\item"c)" if $\beta   -1-\varepsilon
                          \leq \alpha  \leq \beta   -2$, $\alpha  \leq
d-\varepsilon    -3$, then $c(H,g^r_n)\leq \alpha  +1$. \endroster
In
particular, if $dim \ g^r_n(-H)=r(n-d,d)-\varepsilon '$, in the cases a)
with $\alpha  >\beta   -1$ and b) we have $\varepsilon '\leq
\varepsilon$, in the cases a) with $\alpha  =\beta   -1$ and c) we have
$\varepsilon    '\leq \varepsilon    +1$.\endproclaim
\demo{Proof} a)
Let us assume by contradiction that $c(H,g^r_n)\geq     \alpha  +2$;
then, by \cite C,1.3, if $i\leq d- \alpha  - 2$, $dim \ (g^r_n\oplus
\mid iH\mid )-     dim \ (g^r_n\oplus    \mid (i-     1)H\mid )\geq
i+2+\alpha$. In particular, $$\split dim (g^r_n\oplus    \mid (d-
\alpha  -     3)H\mid ) &\geq      \frac{d(d-     1)}2 -
\frac{(\alpha  +2)(\alpha  +3)}2 + r(n,d)-     \varepsilon     = \\
&=  \frac{d(d-     1)}2 -     (\alpha  +3+\varepsilon    )-     \beta
   =\\  &=p_a(C)+d-     (\alpha  +4+\varepsilon    )-     \beta
\geq \\&\geq     p_a(C)-     \beta.\endsplit$$ Let $E$ be a general
divisor of $g^r_n\oplus    \mid (d-     \alpha  -     3)H\mid$; then the
degree of $E$ is $deg \ E=2p_a(C)-     2-     \beta$    and the index of
speciality $i(E)\geq     d-     2-     \alpha  -     \varepsilon
\geq     1$. So, if $d>\alpha  +\varepsilon    +3$, the series $\mid
K-     E\mid$  has degree $\beta    \leq d-     \varepsilon    -     2$
and dimension at least one, which is impossible; let $d=\alpha
+\varepsilon    +3$: if $i(E)=2$ the conclusion is the same, if $i(E)=1$,
$\mid K-     E\mid$  is a fixed divisor $F$ of degree $\beta$, so
$g^r_n\oplus    \mid \varepsilon    H\mid = \mid (d-     3)H-     F\mid$
and $g^r_n= \mid (d-     3-     \varepsilon    )H-     F\mid  = \mid
\alpha  H-     F\mid$; this is impossible because $r = \frac 12 \alpha
(\alpha  +3)-     \beta   -     \varepsilon$.

 b) In
this case, if, by contradiction, $c(H,g^r_n)\geq     \alpha  +1$, arguing
as in a) we find $$\split dim \mid E\mid  &\geq      p_a(C) -
\frac{(\alpha  +1)(\alpha  +2)}2  + \frac{(\alpha  -     1)(\alpha
+2)}2 -      \varepsilon     =\\&= p_a(C) -      (\alpha
+2+\varepsilon    ) \geq\\&\geq      p_a(C)-     \beta,\endsplit$$ so
$dim  \mid K-     E\mid  \geq     1$; since $deg \mid K-     E\mid  =
\beta    \leq  d-     1$, then the only possibility is $dim \mid K-
E\mid  =1$, $\beta   =d-     1$ and $\mid K-     E\mid =\mid H-
P\mid$. This cannot happen if $\beta   >\alpha  +2+\varepsilon$; if
$\beta   =\alpha  +2+\varepsilon$, then $g^r_n=\mid (\alpha  -
1)H+P\mid$, which is impossible because $r=\frac 12 (\alpha  -
1)(\alpha  +2)-     \varepsilon$.

c) If
$c(H,g^r_n)\geq     \alpha  +2$, $\alpha  \leq \beta   -     1$, then
$$\split dim  \mid E\mid  &\geq      \frac{d(d-     1)}2 -
\frac{(\alpha  +2)(\alpha  +3)}2 +\frac{(\alpha  -     1)(\alpha  +2)}2
-     \varepsilon     =\\&= \frac{d(d-     1)}2 -     2(\alpha  +2)-
\varepsilon.\endsplit$$ So, by a direct computation, we get: $dim \mid
E\mid  -     (p_a(C)    -     1-     \beta   ) \geq      d-     4 -
2\alpha  +\beta   -     \varepsilon     \geq      \beta   -     \alpha
-     1$, because $\alpha  \leq d-     \varepsilon    -     3$. Since
$\beta   >\alpha  +1$, then we get $dim \mid K-     E\mid \geq     1$ and
$deg  (K-     E)=\beta$    which is a contradiction.\enddemo
\demo{Proof of
Theorem 3.2} We assume $\varepsilon    =1$ and proceed by induction on
$\alpha$.

The first step of the induction is for $\alpha  =2$. If $\beta
\geq     3$, then $r=1$; let $D$ be a general divisor of the series
$g^1_n$: we have $c(D,K) = dim \mid K\mid -     dim\mid K-     D\mid  =
p_a(C)    -     (1-     n+p_a(C)) = n-    1$, so $D$ imposes $n-     1$
conditions to the curves of degree $d-     3$. If $D$ is the sum of
$n$ distinct points, we apply 3.5 with $\tau =d-     3$: if $\beta
\geq     4$, then we may take $s=2$, $t=1$, so one of the following
happens:

\roster\item"-"  $D$ contains $d-     1$ points on a line, so
$g^1_n    = \mid H-     P+P_1+...+P_{d-     \beta   +1}\mid$  (case ii1))

\item"-"  $\beta
=4$ and $D$ is contained in a conic, so $g^1_n    = \mid 2H-     P_1-
...-     P_4\mid$  (case i2)); \endroster
if $\beta   =3$, then we may take
$s=3$ and there are the following
possibilities:\roster\item"-"     $D$ contains $d-     1$ points
on a line ($g^1_n    = \mid H-     P+P_1+...+P_{d-     2}\mid$, case
ii1)),
\item"-"      $2d-     3$ points of $D$ lie on a conic, so  $g^1_n     =
\mid 2H-     P_1-     ...-     P_4+P\mid$  (case i1)),
\item"-"      $d=6$ and
$D$ is a complete intersection of type $(3,3)$. In this case $$g^1_n
=\mid 3H-     P_1-     ...-     P_9\mid,$$ where $P_1,...,P_9$ impose $8$
conditions to the cubics (we are in the case (iii)). \endroster

If $D$ is not the sum of distinct points, then $g^1_n$     has fixed
points; so we remove them and reduce to one of the previous cases.

If $\beta   =0,1,2$, by lemma 3.7 a), $\varepsilon    '\leq 1$.

If
$\beta   =0$, the series is a $g^4_{2d}$     such that $g^4_{2d}(-
H)$ has dimension $2$ or $1$: the first case is impossible, because it
would imply $g^4_{2d}(-     H)=\mid  H\mid$   and $g^4_{2d}=\mid
2H\mid$; in the second case by 3.6 $g^4_{2d}(-     H)=\mid H-
P+Q\mid$, so $g^4_{2d}=\mid 2H-     P+Q\mid$  (case i1)).

If $\beta
=1$, the series is a $g^3_{2d-1}$  such that $g^3_{2d-1}(-     H)$ either
is a fixed divisor $F$ of degree $d-     1$, or is of the form $\mid
H-     P\mid$. In the first case $g^3_{2d-1}= \mid H+F\mid$, where $F$ is
not fixed; by 3.4, $c(F, (d-     4)H) = d-     2$ so, by 3.5, $F$
contains $d-     2$ points lying on a line, i.e. $F=\mid H-     P-
Q+R\mid$, then $g^3_{2d-1} = \mid 2H-     P-     Q+R\mid$  (case i1)). The
second case is clearly impossible.

If $\beta   =2$ and the series
$g^2_{2d-2}$ is such that $g^2_{2d-2}(-     H)\neq\emptyset$, then
$g^2_{2d-2}= \mid H+P_1+...+P_{d-     2}\mid$  (we are in the case ii2));
if $g^2_{2d-2}(-     H) = \emptyset$, then $c(D,K) = n-     2$ and
$$\split c(D, \mid (d-     2)H\mid ) &= dim J\ \mid (d-     2)H\mid
-      dim \ \mid (d-     2)H-     D\mid  =\\&= \frac{(d-     2)(d+1)}2
-     (p_a(C)    -     d)= n;\endsplit$$ so $h^1(\Cal I_D(d-     3))=2$,
$h^1(\Cal I_D(d-     2))=0$. Let $\chi (D)$ denote the numerical
character of $D$ (see \cite {GP2} for the definition and first
properties); we have: $\chi (D)=(d-     1, d-     1, n_2,...)$, where
$2d-     2=(d-     1)+(d-     2)+(n_2-     2)+...$, so $\chi (D)=(d-
1,d-     1,3)$. Since by assumption $d>5$, the character is disconnected
and, by \cite{JEP}, $2d-     3$ points of $D$ are contained in a conic;
we conclude that   $$g^2_{2d-2}= \mid 2H-     P_1-     P_2-
P_3+Q\mid$$  (case i1)).

Assume now $\alpha  >2$. If $\alpha  >\beta   -     1$, $\alpha
\leq d-     \varepsilon    -     3$, by lemma 3.7 $dim J\ g^r_n(-
H)\geq     r(n-     d,d)-     1$ and $\varepsilon    '=0$ or $1$. By
induction one of the following happens: \roster
\item $g^r_n(-     H) = \mid
(\alpha  -     1)H-     D'+F'\mid$, $deg \ D' = \beta   +\varepsilon
'$, $deg \ F'=\varepsilon    '$;
\item $g^r_n(-     H) = \mid (\alpha  -
2)H-     D'+F'\mid$, $deg \ D' = \beta   +\varepsilon    '-     \alpha$,
$deg \ F' =d+\varepsilon    '-     \alpha$, $\alpha  \leq \beta
+\varepsilon    '$.\endroster
In the case (1),  $g^r_n = \mid \alpha  H-
D'+F'\mid$, where $D'$ imposes independent conditions to the curves of
degree $\alpha  -     1$, hence also to the curves of degree $\alpha$;
since $r=\frac 12 \alpha  (\alpha  +3)-     \beta   -     1 \leq
\frac 12 \alpha  (\alpha  +3) -      deg D'$, we get $\varepsilon
'=1$, $c(H,g^r_n)=\alpha  +1$ and the series is of type i1). In the case
(2), $g^r_n=\mid (\alpha  -     1)H-     D'+F'\mid$; if $F$' is fixed
for $g^r_n$, we get $\frac{(\alpha  -     1)(\alpha  +2)}2 -      deg \
D' = \frac {\alpha  (\alpha  +3)}2 -     \beta   -     1$, so $deg \
D'=\beta   +\varepsilon    '-     \alpha$, $\varepsilon    '=0$ (the
series is of type ii1)     or ii2)). If $F'$ is not fixed,
$\varepsilon    '=1$ and $\beta   \leq \alpha  \leq \beta   +1$. By 3.4,
$F'$ imposes $d-     \alpha  = deg J F'-     1$ conditions to $\mid K-
(\alpha  -     1)H+D'\mid  = \mid (d-     \alpha  -     2)H+D'\mid$, so
by 3.5 $d-     \alpha$   points of $F'$ lie on a line, hence $g^r_n=\mid
\alpha  H-     D+Q\mid$  (case i1)).

If $\beta   \geq     \alpha  +3$, $dim
 \ g^r_n(-     H)\geq     r(n-     d,d)-     1$; by induction
$g^r_n(-     H)= \mid (\alpha  -     2)H-     D'+F'\mid$,
$deg \ D'=\varepsilon    '$, $deg \ F'=d-     \beta   +\varepsilon    '$;
hence $g^r_n= \mid (\alpha  -     1)H-     D'+F'\mid$, which implies that
$\varepsilon    '=1$ and the series is of type ii1).

If $\beta
=\alpha  +1$, $\alpha  +2$, then $\varepsilon    '\leq 2$. If
$\varepsilon    '\leq 1$, we may proceed by induction. If $\beta
=\alpha  +1$, one of the following happens: \roster
\item $g^r_n(-     H) =
\mid (\alpha  -     1)H-     D'\mid$, $deg \ D' = \beta$    and
$\varepsilon    '=1$;
\item $g^r_n(-     H) = \mid (\alpha  -     2)H-
D+F\mid$, $deg \ D' = \varepsilon    '$, $deg \ F' =d-     \beta
+\varepsilon    '$.\endroster
 The first case is impossible because it
implies $g^r_n= \mid \alpha  H-     D'+F'\mid$  where $D'$ is a divisor of
degree $\alpha  +1$ which imposes dependent conditions to the curves of
degree $\alpha$; in the second case $g^r_n= \mid (\alpha  -     1)H-
D'+F'\mid$  and $\varepsilon    '=1$ (case ii1)). If $\beta   =\alpha
+2$, only the case (2) for $g^r_n(-     H)$  is possible; we get
$g^r_n= \mid (\alpha  -     1)H-     D'+F'\mid$, $\varepsilon
'=1$ (case ii1)).

Assume $\varepsilon    '=2$, then $c(H,g^r_n)=\alpha  +1$. With notations
as in lemma 3.7, we have: $dim \mid E\mid  \geq      p_a(C)    -
(\alpha  +3)$, $i(E)\geq      \beta  -     \alpha  - 1$. If $dim \mid
K-     E\mid \geq     1$, then $\beta   =d-     1$, $\mid K-     E\mid  =
\mid H-     P\mid  = \mid \alpha  H-     g^r_n\mid$  so $g^r_n= \mid
(\alpha  -     1)H+P\mid$  which is impossible.

If $dim \mid K-     E\mid
=0$, $\mid K-     E\mid =F$, fixed divisor of degree $\beta$
and\linebreak $g^r_n= \mid \alpha  H-     P_1-     ...-     P_{\beta}
\mid$; since this series is of dimension at least $\frac{\alpha  (\alpha
+3)}2 -     \beta$, it follows $\beta   =\alpha  +2$ (case i2)).

 Finally assume
$\beta   =\alpha  +1$, $\mid K-     E\mid  =\emptyset$, $dim \mid E\mid
= p_a(C)    -     (\alpha  +3)$. Consider $E_1= g^r_n\oplus    \mid
(d-     \alpha  -     4)H\mid$; since $c(H,g^r_n)=\alpha  +1$, $dim \mid
E_1\mid \geq      \frac{(d-     3)(d-     2)}2 -(\alpha  +3)$, moreover
equality holds because it holds in the analogous expression for $E$. We
get: $deg \ E_1= d(d-     4)-     \beta$, $deg \ (K-     E_1)=2d-
(d-     \beta   )$, $i(E_1)=2$ so we may apply the first step of the
induction to $\mid K-     E_1\mid$. Since $d-     \beta   \geq     3$,
$\mid K-     E_1\mid$  has one the following forms:

 - $\mid K-
E_1\mid  = \mid H-     P+P_1+...+P_{\beta   +1}\mid$, or

- $\mid
K-     E_1\mid  = \mid 2H-     P_1-     ...-P_4\mid$  and $d-     \beta
=4$, or

  - $\mid K-     E_1\mid  = \mid 2H-     P_1-     ...-P_4+P\mid$
and $d-     \beta   =3$.

In the first case, $\mid (\alpha  +1)H-
g^r_n\mid  = \mid H-     P+P_1+...+P_{\beta   +1}\mid$, so
\linebreak $g^r_n= \mid \alpha  H-     P_1-     ...-     P_{\beta
+1}+P\mid$  (case i1)). The second case leads to a contradiction, because
we would have $g^r_n=\mid (\alpha  -     1)H+P_1 + ... +P_4\mid$  with
$r=\frac{(\alpha  -     1)(\alpha  +2)}2 -     1$. Finally, in the third
case we get  $g^r_n= \mid (\alpha  -     1)H -     P+ P_1 + ...
+P_4\mid$   (case ii1)).\enddemo

 \head 4. Lifting theorems. \endhead

Let $X$ be an integral non--degenerate variety of codimension $2$ in
$\Bbb P^{r+2}$ and of degree $d$; let $\sigma$ be the minimal degree
of a degenerate hypersurface containing its general hyperplane
section $Y$. From now on we make the assumption:
$d>\sigma ^2 - (r - 1)\sigma +
\binom r2 +1$. We shall prove that, if $r<3$ or $r=3$,
$\sigma >5$, this implies $H^0(\Cal I_X(\sigma ))\ne (0)$.  The method is
the following: suppose by contradiction that  $h^0(\Cal I_X(\sigma
))=0$, then, by the results of \S 2, the total focal linear system
$\delta _T$ on $G_Y$ has dimension at least $2r+n-2$; if we cut $G_Y$ with
$r-1$ general hyperplanes, we get an irreducible plane curve
$G_{\Gamma}$; the restriction of $\delta _T$ to $G_{\Gamma}$ is a linear
 series whose variable part has degree at most
 $\sigma (\sigma +1)-d <
r\sigma -\binom r2-1$ and whose dimension we may control. Thanks to the
results of \S 3 on special linear series on plane curves, we are able to
conclude.

The first application of this method is for $r=1$. In this case
we find once more the Laudal's lemma:
\proclaim{4.1. Theorem} (Laudal's lemma) Let
$C$ be an integral curve of $\Bbb P^3$, let $\sigma$ be the minimal
degree of a plane curve containing a general plane section $\Gamma$ of
$C$. If $d>\sigma ^2+1$, $\sigma >1$, then $C$ is contained  in a
surface of degree $\sigma$. \endproclaim
\demo{Proof} By the assumption on $d$,
$n=1$. Let $\delta$ be the focal linear system on $G_{\Gamma}$, a plane
curve of degree $\sigma$ containing $\Gamma$;  its variable part is a
linear series of degree $\sigma (\sigma +1)-d\leq \sigma -2$. If
$H^0(\Cal I_C(\sigma))=(0)$, then by theorem 2.3 $dim \ \delta \geq 1$:
this contradicts 3.1. \enddemo
\remark{4.2. Remark} Let us explicitly note
that, to prove theorem 4.1, it is not necessary to assume that $C$ is
integral; in fact it suffices to assume that $G_{\Gamma}$ is
integral.\endremark
  We will study now the cases with $r\geq2$.
\proclaim{4.3. Lemma}
Let $X$ be an integral variety of dimension $r$ and degree $d$ in
$\Bbb P^{r+2}$, $\Gamma$ its general plane section. If $h^0(\Cal I_{\Gamma
}(\sigma -1))=0$ and
$d \geq \sigma ^2-(r-1)\sigma +\binom r2 +1$, then $h^0(\Cal I_{\Gamma}
(\sigma)) \leq r$.\endproclaim
\demo{Proof} Let $h_\Gamma (t)$ be the Hilbert function of $\Gamma$ and
$\Delta h_\Gamma (t)=h_\Gamma (t)-h_\Gamma (t-1)$ be its first
variation. Assume that $h^0(\Cal I_{\Gamma}(\sigma))\geq r+1$: this means
that $\Delta h_\Gamma (t)=t+1$ for $t<\sigma$ and $\Delta h_\Gamma
(\sigma )\leq \sigma -r$; since $\Delta h_\Gamma (t)$ is strictly
decreasing for $t\geq \sigma$ by the uniform position property (\cite
{Ha}), the assumption on $d$ gives a contradiction. \enddemo

\proclaim{4.4.
Theorem} Let $S$ be an integral non-degenerate surface of $\Bbb P^4$ of
degree $d$. Let $\sigma$ be the minimal degree of a degenerate
surface containing a general hyperplane section of $S$. If $d>\sigma
^2-\sigma +2$, then $S$ is contained in a hypersurface of degree $\sigma$.
\endproclaim
\demo{Proof} Assume first $\sigma\geq 3$. Let $C=S\cap H$ be a general
hyperplane section of $S$ and $\Gamma =C\cap H'$ be a general plane
section of $C$; note that, by 4.1, $h^0(\Cal I_{\Gamma}(\sigma -1))=0$
because $d>(\sigma -1)^2+1$, so by lemma 4.3 $n\leq 2$. Let $\delta _T$
be  the
focal system of the total family on $G_C$, a surface of degree $\sigma$
in $H$ containing $C$, and $\delta _\Gamma =\delta _T\mid_{G_{\Gamma }}$
its restriction to $G_{\Gamma}=G_C\cap H'$. If  $h^0(\Cal I_S
(\sigma ))=0$, then by 2.3 $dim \ \delta _T\geq 2+n$. Let us
consider the exact sequence: $$0 \rightarrow   H^0(\Cal
I_{C,G_C}(\sigma)) @> \cdot H' >> H^0\Cal I_{C,G_C}(\sigma +1))
@> \psi >> H^0(\Cal I_{\Gamma ,G_{\Gamma}}(\sigma +1)).$$
 If $n=1$ $\psi$ is
injective, so $dim \delta _\Gamma = dim \ \delta _T  \geq  3$; if $n=2$,
$dim \ ker\psi\leq 1$, so $dim \ \delta _\Gamma \geq dim \ \delta _T  -1
\geq 3$. But  $deg \ \delta _\Gamma <2\sigma -2$, because $d>\sigma
 ^2-\sigma +2$, so by 3.1 $dim \ \delta _\Gamma \leq 2$: a contradiction.

If $\sigma=2$, $\sigma ^2-\sigma +2=\sigma ^2$, so $n=1$; if in the
homogenous ideal $I_C$ of $C$ there is no cubic, not multiple of $G_C$,
then the focal system on $G_C$ is undetermined, so by 1.8 $G_C$ may be
lifted to a generator of the ideal of $S$: this certainly happens if
$d>6$, by Bezout's theorem. Let $I_C$ contain a new cubic: if $d=6$, then
$C$ is a complete intersection, so it is arithmetically Cohen-Macaulay and
also $S$ is; if $d=5$, then $C$ is linked to a line in a complete
intersection $(2,3)$, so again $C$ and $S$ are arithmetically
Cohen-Macaulay. \enddemo

\remark{4.5. Remarks}  1) In \cite{MR} the above theorem is proved under
the assumption \lq\lq $S$ linearly normal''.  Note that, if $n=2$, this
assumption is not used in the proof and that it is moreover proved that
both generators of $H^0(\Cal I_C(\sigma ))$ are lifted to $H^0(\Cal I_S
(\sigma ))$.

2)  Also in this case we note that, to prove the theorem, it
is not necessary to assume that $S$ is integral; in fact it suffices to
assume that $G_{\Gamma}$, $G_C$ are integral and $n\leq 2$.\endremark

Let us suppose now $r=3$. Let $X$ be an integral non-degenerate $3$-fold
of $\Bbb P^5$ of degree $d$. In the following we  let $S=X\cap H$ be a
general hyperplane section of $X$, $C=S\cap H'$ be a general hyperplane
section of $C$ and $\Gamma =C\cap H''$ be a general plane section of
$C$. As usual, we denote by $\sigma$ the minimal degree of a
hypersurface in $H$ containing $S$.
\proclaim{4.6. Lemma} If $d > \sigma
^2-2\sigma +4$, then $h^0(\Cal I_{\Gamma}(\sigma -1)) = 0$.\endproclaim

\demo{Proof} Otherwise, by 4.1 and 4.4, $S$ would be contained in a
hypersurface of degree $\sigma-1$, because both $d > (\sigma -1)^2+1$ and
$d > (\sigma -1)^2-(\sigma -1)+2$. \enddemo

\proclaim{4.7. Lemma} If $d > \sigma
^2-2\sigma +4$, then $h^0(\Cal I_{\Gamma}(\sigma)) \leq  3$; if $h^0(\Cal
I_{\Gamma}(\sigma)) =  3$, then $h^0(\Cal I_{\Gamma}(\sigma +1)) = 8$. In
particular, if $G_{\Gamma}$ is a curve of degree $\sigma$ containing
$\Gamma$, the elements of $h^0(\Cal I_{\Gamma}(\sigma +1))$ cut on
$G_{\Gamma}$ a linear system of dimension $4$.\endproclaim
\demo{Proof} The
first assertion follows from 4.3 and 4.6. As for the second one, note
that if $h^0(\Cal I_{\Gamma}(\sigma ))=3$ then $\Delta
h_{\Gamma }(\sigma )=\sigma -2$. Since $\Delta h_\Gamma$ is strictly
decreasing from $\sigma$ on and $d>\sigma ^2-2\sigma +4$, we conclude that
$\Delta h_\Gamma (\sigma +1)=\sigma -3$, which gives the lemma. \enddemo
\proclaim{4.8. Lemma} Let $X$ be an integral non--degenerate $3$-fold of
$\Bbb P^5$ of degree $d$ with $\sigma >5$. With notations as above,
let $G_S$ be a general $3$-fold of degree $\sigma$ in $H$ containing $S$,
$\delta$ be the focal linear system of the total family on $G_S$. Let
$\delta _\Gamma$ be the  linear series on $G_{\Gamma}= G_S\cap H'\cap
H''$ obtained by removing the points of $\Gamma$ from $\delta\mid
_{G_\Gamma}$.

If $d>\sigma ^2-2\sigma +4$, then $dim \ \delta _\Gamma \leq 4$; if
equality holds, then one of the following happens: \roster
\item"(i)" $h^0(\Cal
I_{\Gamma}(\sigma )) = h^0(\Cal I_C(\sigma )) = h^0(\Cal I_S (\sigma))
=3$; or
\item"(ii)" $h^0(\Cal I_C(\sigma)) = h^0(\Cal I_S (\sigma)) =2$;
or
\item"(iii)" $h^0(\Cal I_C(\sigma )) = h^0(\Cal I_S (\sigma ))
=1$.\endroster
 In cases (ii) and (iii) moreover $d=\sigma ^2-2\sigma +5$.\endproclaim
\demo{Proof} Since $deg \ \delta _\Gamma = \sigma(\sigma +1)-d <
3\sigma -4$, then, by 3.1, $dim \ \delta _\Gamma \leq 5$.

 If $dim \ \delta _\Gamma =5$,
then the series is complete and, by 3.2, we have $\mid   (\sigma
+1)H-\Gamma\mid    = \mid 2H+F\mid$, where $F$ is a fixed divisor; this
implies that $\Gamma +F\in\mid   (\sigma -1)H\mid$, which contradicts
4.6.

If $dim \ \delta _\Gamma    =4$, by 3.2 one of the following happens:
\roster
\item"a)"
$\mid   (\sigma +1)H-\Gamma\mid    = \mid   2H-P+F\mid$, where $P$ is a
point and $F$ a fixed divisor with $0 \leq  deg \ F \leq \sigma -4$;
\item"b)" $\mid
(\sigma +1)H-\Gamma\mid   =\mid   3H-P_1-\dots -P_5\mid$    where
$P_1,\dots, P_5$ are five points not all on a line and $deg \ \delta
_\Gamma     =3\sigma -5$.\endroster
 In the case a), $\Gamma +F$ is linked in a complete intersection
$(\sigma,\sigma +1)$ to a scheme $\Gamma '$, of degree $2\sigma -1$, which
varies in a linear system of dimension $4$ and is linked in a complete
intersection $(2,\sigma )$ to the point $P$.  By mapping cone
(\cite{PS}, Proposition 2.5), from the minimal free resolution of the
ideal sheaf of $P$ $\Cal I_P$, we get the following free resolutions of
$\Cal I_{\Gamma '}$ and $\Cal I_{\Gamma +F}$ (where $\Cal O$ denotes the
structural sheaf of the plane of $\Gamma$):
$$\split 0& \rightarrow   \Cal O(-2) \rightarrow   2\Cal O(-1)
\rightarrow  \Cal I_P\rightarrow  0,\\
 0& \rightarrow   2\Cal O
(-\sigma -1) \rightarrow   2\Cal O (-\sigma )\oplus     \Cal O(-2)
\rightarrow  \Cal I_{\Gamma '} \rightarrow  0,\\
 0& \rightarrow   \Cal O
(-2\sigma +1)\oplus     2\Cal O    (-\sigma -1)\rightarrow   \Cal O
(-\sigma -1)\oplus     3\Cal O(-\sigma) \rightarrow  \Cal
I_{\Gamma +F} \rightarrow  0.\endsplit
$$
 The last one can be simplified (\cite{PS}, \S 3), to
get the minimal free resolution:
$$0 \rightarrow   \Cal O
(-2\sigma +1)\oplus     \Cal O    (-\sigma -1)\rightarrow   3\Cal O
(-\sigma ) \rightarrow  \Cal I_{\Gamma +F} \rightarrow  0;$$ in
particular  $h^0(\Cal I_{\Gamma}(\sigma )) = 3$. Let us consider now the
curve $\tilde      C$ linked to $C$ in a general complete intersection
$(\sigma, \sigma +1)$: its general plane section is $F+\Gamma '$; since
$F$ is fixed as the liaison varies, there is a fixed curve $\Cal F$
contained in $\tilde      C$ whose general plane section is $F$. Let now
$C'$ be the curve linked to $C+\Cal F$ in a general complete intersection
$(\sigma, \sigma +1)$: its plane section is $\Gamma '$; $\Gamma '$ is
contained in a conic that we may assume to be irreducible and $deg \
\Gamma '>5$, so by 4.2 $C'$ is contained in a quadric; hence $C'$ is
linked to a line in a complete intersection of type $(2,\sigma )$. Arguing
as before, by mapping cone we find that $C$ is contained in $3$ surfaces
of degree $\sigma$. Finally we may repeat the same argument of liaison
starting from $S$ and conclude that also $h^0(\Cal I_S (\sigma )) = 3$.

 In
the case b), $\Gamma$ is linked in a complete intersection
$(\sigma,\sigma +1)$ to a scheme $\Gamma '$, which  is linked in a
complete intersection $(3,\sigma )$ to $\Gamma ''=\lbrace P_1,\dots
, P_5\rbrace$. We may assume that $\Gamma ''$      is contained in an
irreducible cubic: otherwise, by Bertini theorems, $\Gamma ''$
contains $4$ points on a line $l$, so the cubics of the linear system
$\mid   3H-P_1-\dots -P_5\mid$    split in the union of the line $l$ and a
conic through the remaining point and we are again in case a). Let $C'$ be
a curve linked to $C$ in a general complete intersection
$(\sigma,\sigma +1)$; since $deg \ \Gamma '=3\sigma -5>10$ by the
irreducibility assumption, $C'$ too is contained in an irreducible cubic:
let $C''$ be a curve of degree $5$ linked to $C'$ in a complete
intersection $(3,\sigma )$. In an analogous way, the surface $S'$ linked
to $S$ in a general complete intersection $(\sigma,\sigma +1)$ lies on an
irreducible cubic, so it is linked to a quintic surface $S''$ in a
complete intersection $(3,\sigma )$.

Assume first that $C''$ does not lie on a quadric; then
$\Cal I_{C''}$ has a locally free resolution of the form:
$$ 0
\rightarrow   \Cal E \rightarrow   \bigoplus_{i}\Cal O(-n_i)
\rightarrow  \Cal I_{C''}\rightarrow  0, \qquad n_i\geq 3.\tag 4.9 $$
 By applying the
mapping cone twice, we find a resolution of $\Cal I_{C''}$ which ends as
follows:
$$\dots\rightarrow   \oplus_{i}\Cal O
(-n_i-\sigma +2)\oplus     \Cal O    (-\sigma -1)\oplus     \Cal O
(-\sigma ) \rightarrow  \Cal I_C\rightarrow  0,$$
  so $h^0(\Cal I_
C(\sigma ))=1$.

Assume now that $C''$ lies on a quadric. If
it lies on an irreducible quadric, then $S''$ too lies on
a quadric by theorem 4.4 and remark 4.5. If $C''$ lies on a reducible
quadric, then it splits in the union of a plane conic $C_1''$ and a plane
cubic $C_2''$: hence $S''$, which has $C''$ as general hyperplane section,
has to split in the union of two degenerate surfaces, so it lies on a
quadric. In both cases, in the resolution (4.9) and in the analogous one
for $\Cal I_S''$  $n_i=2$ for some $i$, so, again by mapping cone, we
find  $h^0(\Cal I_C(\sigma )) = h^0(\Cal I_S (\sigma )) =2$.\enddemo

\proclaim{4.10. Theorem}  Let $X$ be an integral non--degenerate
subvariety of $\Bbb P^5$ of dimension $3$ and degree $d$. Let $\sigma$
be the minimal degree of a degenerate hypersurface containing a general
hyperplane section of $X$. If $d>\sigma ^2-2\sigma +4$ and $\sigma >5$,
then $X$ is contained in a hypersurface of degree $\sigma$. Moreover, if
$h^0(\Cal I_C(\sigma )) = h^0(\Cal I_S (\sigma ))=2$, both generators lift
to $H^0(\Cal I_X(\sigma ))$, if $h^0(\Cal I_S (\sigma))=3$, all the three
generators lift to $H^0(\Cal I_X(\sigma ))$.\endproclaim
\demo{Proof} Assume
first that $h^0(\Cal I_C(\sigma )) = h^0(\Cal I_S (\sigma )) =2$. Fixed a
general line $l$, we construct the family $\Sigma _Z$, $Z=\check\Bbb
P^5$, as in \S 1. Let $H$ have equation $x_5=0$, and $l\cap
H=P(0:\dots :0:1:0)$; then a focal matrix of $H$ in $\Sigma _Z$ has the
form: $M=\left( \matrix  F_0  & F_1 &  \hdots & F_4\\
 x_0 & x_1 & \hdots & x_4  \endmatrix \right).$ Let  $\Sigma _{Z'}$
be the subfamily of $\Sigma _{Z}$  parametrized by the set of
 all hyperplanes through $P$: a focal matrix of $H$ in
$\Sigma _{Z'}$ is $N= \left( \matrix  F_0   &  \hdots & F_3\\
 x_0  & \hdots & x_3  \endmatrix \right);$ the minors of $N$ define a
subsystem $\delta '$ of the focal linear system $\delta$ associated to
$\Sigma _{Z}$      on $G_S$, which is the image via the map $\phi _M$,
introduced in \S 2,  of the planes through $P$; these planes form a
Schubert cycle $\Lambda$ of codimension $2$ in the grassmannian $\Bbb
G(2,4)$ which has dimension $6$. Assume by contradiction that there is no
hypersurface of $\Bbb P^5$ of degree $\sigma$ containing $X$ and passing
through $P$; then, we may argue as in the proof   of
Theorem 2.3: by considering the intersection of $\Lambda$ with a linear
space of codimension $3$, we conclude that $dim \ \delta '\geq 3$. Fix
now $A\in H^0(\Cal I_S (\sigma ))$, $A\neq G_S$ and consider the linear
system $\delta ''$ cut out on $G_S$ by the hypersurfaces of degree
$\sigma +1$ through $S$; it contains the multiples of $A$ and $6
$ independent minors of $M$, $4$ of them extracted from $N$. Note that the
polynomials $F_i$, $i=0,\dots ,3$,  vanish at $P$, because they correspond
to infinitesimal deformations fixing $P$, so the maximal minors of $N$
have a zero of order at least two at $P$, modulo $G_S$. On the other hand
only the product of $A$ with the tangent hyperplane to $G_S$ at $P$
vanishes doubly at $P$, so at least $4$ multiples of $A$ are independent
of the minors of $N$ and $dim \ \delta ''\geq7$; by cutting twice with
hyperplanes, we find a linear system on $G_{\Gamma}=G_S\cap H\cap H'$ of
dimension $5$ and degree $< 3\sigma -4$, which is impossible by 4.8. So
we have proved that there is a hypersurface of degree $\sigma$
containing $X$ and passing through $P$.

 Let us assume now that $n=3$: a
focal matrix of $H$ associated to the total family $\Sigma _T$ has the
form: $M = \left( \matrix  F_0  & F_1 &  \hdots & F_4 & A & B\\
 x_0 & x_1 & \hdots & x_4 & 0 & 0 \endmatrix \right),$ where
$A,B\in H^0(\Cal I_S (\sigma ))$ by 1.8. The multiples of $A$ and $B$
generate on $G_S$ in degree $\sigma +1$ a linear system $\delta '$ of
dimension at least $8$; if the dimension were $9$, by cutting successively
with hyperplanes we would find a linear system of degree
$\sigma(\sigma +1)-d<3\sigma -4$ and dimension at least $5$ on
$G_{\Gamma}$, which is impossible by 4.8. So $dim \ \delta '=8$. Let us
fix two general points $P, Q$ on $G_S$: we may assume $P(0:1:0:\dots:0)$
and $Q(1:0:\dots:0)$; consider the subfamily $\Sigma '$ of $\Sigma _T$
formed by the hypersurfaces lying on hyperplanes containing $P$ and $Q$
and passing through $P$ and $Q$. A focal matrix for $\Sigma '$ has the
form $N=\left( \matrix  F_2  & F_3  & F_4\\
 x_2 & x_3  & x_4  \endmatrix \right).$ As above, note that $F_i$,
$i=2,3,4$, vanish at $P$ and $Q$, because they correspond to infinitesimal
deformations fixing $P$ and $Q$, so the maximal minors of $N$ have a zero
of order at least two at $P$ and $Q$, modulo $G_S$. Moreover these minors
are in $\delta '$,  because otherwise the dimension of the focal linear
system would be $9$. Let us prove that   only one element of $\delta '$
may have a zero of order two both in $P$ and in $Q$: in fact, a double
point in fixed position imposes $4$ conditions to a general element of
$\delta '$; if for any pair of points $(P,Q)$ the $8$ conditions of having
double points at $P$ and $Q$ were dependent, then, letting $(P,Q)$ move on
$G_S$, we would have a codimension one subvariety $V$ of $\delta '$; the
tangent space to $V$ at a point corresponding to a pair of distinct points
$(P_0,Q_0)$ parametrizes elements of $\delta '$ passing through  $P_0$ and
$Q_0$, so it should have codimension $2$ in $\delta '$: a contradiction.
So we have proved that the minors of order two of $N$ are two by two
proportional; now we show that they are zero (mod $G_S$). Assume for
example that $F_2x_3-F_3x_2=a(F_2x_4-F_4x_2) (mod G_S)$, $a\in JK$; then
$F_2(x_3-ax_4)=x_2(F_3-aF_4)(mod G_S)$, which implies that $F_2$ is a
multiple of $x_2$, mod $G_S$. We get that $N$ has the form $N=\left(
\matrix  A_2x_2  & A_3x_3 &  A_4x_4\\
 x_2 & x_3 &  x_4  \endmatrix \right).$ The minors $(A_i-A_j)x_ix_j$,
$i,j=2,3,4$, $i\neq j$, are zero along $S$; by the minimality of
$\sigma$, it follows that $A_i=A_j$. So we conclude that $G_S$ may be
lifted to a hypersurface containing $X$ and passing through $P$ and $Q$.

Let us finally assume that $h^0(\Cal I_X(\sigma ))=0$, $n=1,2$; then
$dim \ \delta _T  \geq 4+n$. Consider the exact sequences:
$$0 \rightarrow   H^0(\Cal I_{S,    G_S}(\sigma )) @> \cdot H' >>
H^0(\Cal I_{S,G_S}(\sigma +1)) @> \psi >>  H^0(\Cal
I_{C,G_C}(\sigma +1));$$
$$  0 \rightarrow   H^0(\Cal I_{C,G_C} (\sigma))
@> \cdot H'' >> H^0(\Cal I_{C,G_C}(\sigma +1)) @> \phi >> H^0(\Cal
I_{\Gamma ,G_\Gamma }(\sigma +1)),$$
 where $G_S\in H^0(\Cal I_S (\sigma ))$,
$G_C=G_S\cap H'$ and $G_{\Gamma }=G_C\cap H''$.

 Assume $n=1$: by the above
sequences, if $h^0(\Cal I_C(\sigma ))=1$, then $dim \ \delta _\Gamma
\geq 5$; if $h^0(\Cal I_C(\sigma ))=2$,  then $dim \ \delta _\Gamma
\geq 4$; if $h^0(\Cal I_C(\sigma ))=3$, by 4.7 $$dim \ \Bbb P(h^0(\Cal
I_{\Gamma}(\sigma +1)))\mid  _{G_{\Gamma}}=4;$$ by Lemma 4.8 all these
possibilities are excluded. If $n=2$, the case $h^0(\Cal I_C(\sigma ))=3
$ is excluded as before; if $h^0(\Cal I_C(\sigma ))=2$, $dim \ \delta _T
\geq 6$, so we should have $dim \ \delta _\Gamma    \geq 4$, which implies
$h^0(\Cal I_S (\sigma ))=2$ by 4.8. In this case we have proved that
$h^0(\Cal I_X(\sigma ))=2$.\enddemo


\Refs
\widestnumber\key{GP2}



\ref\key Ch\by	M. C. Chang \paper Characterization of arithmetically
Buchsbaum subschemes of codimension $2$ in $\Bbb P^n$ \jour J. Diff. Geom.
\vol 31 \yr 1990 \pages 323--341  \endref
\ref\key
C\by C. Ciliberto \paper Alcune applicazioni di un classico procedimento
di Castelnuovo  \inbook Seminari di Geometria \publ Universit\`a di
Bologna \yr 1982--83 \pages 17--43 \endref
\ref\key CS\by	C. Ciliberto, E. Sernesi \paper Singularities of the
theta divisor and congruences of planes \jour  J. Algebraic
Geometry \vol 1  \yr 1992
\pages 231--250  \endref
\ref\key E\by	Ph. Ellia
\paper Sur le genre maximal des courbes gauches de degr\'e $d$ non sur une
surface de degr\'e $s-1$ \jour J. reine angew. Math. \vol 413 \yr 1991
\pages 78--87   \endref
\ref\key EP\by
Ph. Ellia, Ch. Peskine \paper Groupes de points de $\Bbb P^2$:
caract\`ere et position uniforme \inbook Lecture Notes in Math. \vol 1417
\publ   Springer-Verlag
\publaddr Berlin, Heidelberg, and New York
\yr     1990
\pages  111--116
\endref
\ref\key
GP\by L. Gruson, Ch. Peskine \paper Section plane d'une courbe gauche:
postulation \inbook Enumerative Geometry \bookinfo Progr. in Math. \vol
24 \publ Birkh\"auser \yr 1982 \pages 33--35 \endref
\ref\key GP2\by L. Gruson, Ch. Peskine \paper Genre des courbes de
l'espace projectif \inbook Lecture Notes in Math. \vol 687
\publ   Springer--Verlag
\publaddr Berlin, Heidelberg, and New York
\yr     1978
\pages  31--60
\endref
\ref\key H\by	R. Hartshorne \paper Generalized divisors on Gorenstein
curves and a theorem of Noether \jour J. Math. Kyoto Univ. \vol 26--3
\yr 1986 \pages 375--386 \endref
\ref\key L\by O. A. Laudal \paper A generalized trisecant lemma
\inbook Lecture Notes in Math. \vol 687
\publ   Springer-Verlag
\publaddr Berlin, Heidelberg, and New York
\yr     1978
\pages  112--149
\endref
\ref\key M\by	E.
Mezzetti \paper The border cases of the lifting theorem for surfaces in
$\Bbb P^4$ \jour J. reine angew. Math. \toappear \endref
\ref\key MR\by	E. Mezzetti, I. Raspanti
\paper A Laudal--type theorem for surfaces in $\Bbb P^4$ \jour Rend. Sem.
Mat. Univ. Pol. Torino \toappear  \endref
\ref\key N\by	M. Noether \paper Zur Grundlegung
der Theorie der algebraischen Raumcurven \inbook K\"oniglichen Akad. der
Wissenschaften \publaddr Berlin \yr 1883 \endref
\ref\key PS\by	Ch. Peskine, L. Szpiro \paper Liaison des vari\'et\'es
alg\'ebriques \jour Inventiones Math. \vol 26 \yr 1974 \pages 271--302
\endref
  \ref\key Se\by	C. Segre \paper Preliminari di una teoria delle
variet\`a luoghi di spazi\jour Rend. Circ. Mat. Palermo \vol 30 \yr  1910
\pages 87--121 \endref
\ref\key St\by R.
Strano \paper On generalized Laudal's lemma \inbook Projective
Complex Geometry \publ Cambridge University Press \toappear \endref


\endRefs
\enddocument